\documentclass{article}

\usepackage{amsmath,amssymb,amsthm,bm}
\usepackage[a4paper]{geometry}
\usepackage{algorithm,algorithmic}
\usepackage{caption}
\usepackage{graphicx}
\usepackage[colorlinks=false,pdfborder={0 0 0}]{hyperref}
\usepackage{xcolor}
\usepackage[T1]{fontenc}
\usepackage{multirow}
\usepackage{nicematrix}
\usepackage{authblk}
\DeclareMathOperator*{\argmin}{\arg\min}

\DeclareMathOperator{\diag}{diag}
\DeclareMathOperator{\Range}{Range}

\DeclareMathOperator{\Span}{span}

\DeclareMathOperator{\VECz}{vec_0}

\newcommand*{\set}[1]{\left\lbrace#1\right\rbrace}

\newcommand*{\herm}{^*}
\newcommand{\bmat}[1]{\begin{bmatrix}#1\end{bmatrix}}

\newcommand*{\md}{\mathop{}\mathopen{}\mathrm{d}}
\newcommand*{\mi}{\mathrm i}
\newcommand*{\nep}{n_{\rm{ep}}}
\newcommand*{\save}{{\rm{log}}}

\graphicspath{{fig/}}

\newenvironment{keywords}{\medskip\textbf{Keywords:}}{}
\newenvironment{AMS}{\medskip\textbf{AMS subject classifications (2020).}}{}

\newtheorem{theorem}{Theorem}

\newtheorem{lemma}[theorem]{Lemma}

\newtheorem{remark}{Remark}

\definecolor{bkgndcolor}{rgb}{1,1,1}

\title{Improving performance of contour integral-based nonlinear eigensolvers with infinite GMRES}
\author[1]{Yuqi Liu}
\author[2]{Jose E. Roman}
\author[3,4]{Meiyue Shao}
\affil[1]{School of Mathematical Sciences, Fudan University, Shanghai 200433,
China}
\affil[2]{D. Sistemes Inform\`atics i Computaci\'o, Universitat Polit\`ecnica de Val\`encia, Cam\'i de Vera s/n, 46022 Val\`encia, Spain}
\affil[3]{School of Data Science, Fudan University, Shanghai 200433, China}
\affil[4]{Shanghai Key Laboratory for Contemporary Applied Mathematics, Fudan
University, Shanghai 200433, China}
\date{\today}

\begin{document}
\pagecolor{bkgndcolor}

\maketitle

\begin{abstract}
In this work, the infinite GMRES algorithm, recently proposed by Correnty
\emph{et al.}, is employed in contour integral-based nonlinear eigensolvers,
avoiding the computation of costly factorizations at each quadrature node to
solve the linear systems efficiently.
Several techniques are applied to make the infinite GMRES memory-friendly,
computationally efficient, and numerically stable in practice.
More specifically, we analyze the relationship between polynomial eigenvalue
problems and their scaled linearizations, and provide a novel weighting
strategy which can significantly accelerate the convergence of infinite GMRES
in this particular context.
We also adopt the technique of TOAR to infinite GMRES to reduce the memory
footprint.
Theoretical analysis and numerical experiments are provided to illustrate the
efficiency of the proposed algorithm.
\end{abstract}



\section{Introduction}

By nonlinear eigenvalue problem (NEP), we refer to
\begin{equation}
\label{eq:nep}
T(\lambda)v=0,
\qquad v\in\mathbb{C}^n\setminus\set{0},
\qquad \lambda\in\Omega,
\end{equation}
where \(\Omega\subseteq\mathbb{C}\) is a connected region with a smooth
boundary, and \(T(\xi)\colon\Omega\rightarrow\mathbb{C}^{n\times n}\)
is a \(\xi\)-dependent matrix~\cite{GT2017,MV2004}.
In the particular case that \(T(\xi)\) is a matrix polynomial in \(\xi\),
\eqref{eq:nep} is also called a polynomial eigenvalue problem (PEP).
Many well-known problems are of this type, such as the Orr--Sommerfeld
equation~\cite{TH2001} for the incompressible Navier--Stokes equation,
and the vibration analysis of the damped beams~\cite{HMT2008}.

Nowadays, an increasing number of non-polynomial NEPs arise form physics,
chemistry, and industrial applications.
These problems have remarkable practical and theoretical value.
But they are far more complicated to solve due to their large scale and
nonlinearity.
One emerging example appears in quasi normal mode (QNM)
analysis~\cite{NDZ2023}, where the computation of a small set of eigentriplets
of large, sparse rational eigenvalue problems
\[
\sum_{j=1}^pr_j(\lambda)L_jv=0,
\qquad r_j(\cdot)\colon\mathbb{C}\rightarrow\mathbb{C},
\qquad L_j\in\mathbb{C}^{n\times n},
\qquad j=1,\dotsc,p
\]
is needed for modal expansions, whose results can
be used for the physical understanding of nanophotonic devices.
Actually, solving rational eigenvalue problems for QNM analysis is currently a
very active field; see also~\cite{BBC2024,LYG2019,Ming2023}.
Examples with other nonlinear functions can be found in~\cite{HR2023},
where nonlinear eigenvalue problems of the form
\begin{equation}
\label{eq:NEP-ABCS}
(A+\lambda B+\lambda^2C)v=e^{\mi\lambda\tau}Sv
\end{equation}
are to be solved.
Here, \(A\), \(B\), \(C\), and \(S\) are constant matrices, and \(\tau\) is a
scalar.
The solution of~\eqref{eq:NEP-ABCS} can be used to obtain the acoustic modes
of the Helmholtz wave equation with high efficiency and accuracy.

To address these problems,
it may be beneficial to employ
Krylov-based methods.
This kind of algorithms reformulate the NEP as a generalized eigenvalue
problem, and subsequently apply a Krylov algorithm to solve it.
Depending on the specific approximations of \(T(\xi)\), there are Taylor
expansion-based algorithms~\cite{JMM2012}, Chebyshev interpolation-based
algorithms~\cite{KR2014} and, more complicated, rational approximation-based
algorithms~\cite{CR2021,GVM2014,VMM2015}.

Nevertheless, the Krylov-based methods suffer from the drawback that the
eigenvalues have to be distributed in certain patterns.
For example, Taylor expansion proves unsatisfactory when some eigenvalues are
located far away from the expansion point.
Chebyshev interpolation, on the other hand, is specifically designed for the
eigenvalues lying exactly on the real axis or some pre-specified
curves~\cite{EK2012}.
In some practical applications, users seek for the eigenvalues lying in
certain regions, without prior knowledge of their number and distribution.
In such cases, to capture all the eigenvalues by a single
Taylor expansion point is almost impossible, not to say connecting them with
a pre-specified curve.
To address this challenge, researchers
may employ the rational Krylov method
to target multiple points instead of just one.
However, there is still no guarantee that all eigenvalues lying in the region
of interest can be obtained.

Contour integral-based algorithms are another type of widely-used
eigensolvers that are originally developed for linear eigenvalue problems.
One of the most representative algorithms,
the Sakurai--Sugiura method~\cite{SS2003},
transforms the original problems into a generalized eigenvalue problem of two
small Hankel matrices from which it extracts eigenvalue approximations using
standard dense eigensolvers.
Nevertheless, computing these eigenvalues can be numerically
unstable because the Hankel matrices, formed by higher moments, are usually
ill-conditioned.
Thus, projection methods such as FEAST~\cite{Polizzi2009,TP2014} and
CIRR~\cite{SS2007} are developed.
Similar to the Sakurai--Sugiura method,
these methods approximate the characteristic subspace by integration
rules.
However, they subsequently apply a
projection on this subspace to avoid generating Hankel matrices.

When dealing with non-linear eigenvalue problems, the previously mentioned
algorithms have corresponding extensions.
In the case of the Sakurai--Sugiura method, all its theories follow and can
be directly adapted to NEP cases~\cite{AST2009}.
However, the same ease does not apply to FEAST or CIRR, because computing an
approximate solution of a general NEP on a projection subspace is no longer an
easy task.
Hence, nonlinear FEAST~\cite{GMP2018} or CIRR~\cite{YS2013} need additionally
an auxiliary solver for solving projected problems.

To avoid the numerical instability of the Sakurai--Sugiura method
and the uncertainty of nonlinear FEAST and CIRR, we will take Beyn's
algorithm~\cite{Beyn2012} as the eigensolver for this work.
It can be regarded as a nonlinear extension of the Sakurai--Sugiura algorithm
using only the first two moments.
However, it should be noted that, though Beyn's algorithm is employed as
a framework in this paper, the technique we shall propose is not dependent on
a specific eigensolver, but can be extended naturally to any contour
integral-based algorithm where several moments need to be computed.
A more detailed introduction will be provided later.

Although eigenvalues close to the contour may present certain challenges for
contour integral-based algorithms, when compared to Krylov subspace-based
algorithms, contour integral-based algorithms are not significantly reliant on
the a priori knowledge of the eigenvalue distribution.
Hence they are well-suited for extracting all eigenvalues lying in the domain
of interest.
However, these algorithms require solving a series of linear systems, which
could be a heavy workload, especially for large, sparse problems.

Our interest in this work is to solve these linear parameterized systems
efficiently.
If \(T\) is linear with respect to \(\xi\), nested Krylov methods, such as
multi-shift GMRES~\cite{FG1998} or multi-shift QMRIDR~\cite{BV2015} can be
applied to reuse the Arnoldi basis~\cite{AST2010}.
Nevertheless, these techniques are designed for linear cases, and cannot be
easily extended when the system is not linear with respect to the parameter.
In this paper, we employ infinite GMRES~\cite{CJS2023,JC2022} to overcome
these difficulties.
In simple terms, infinite GMRES uses a companion linearization to transform
the parameterized system to a form that is linear with respect
to the parameter, and then, use multi-shift GMRES to solve multiple systems together.

Briefly, our algorithm solves the NEP using contour integral-based algorithms
where the Arnoldi process is employed to solve the linear systems efficiently.
It not only avoids the drawback of poor parallelism in Krylov-based algorithms
by achieving partial parallelism in matrix factorizations and solving least
squares problems, but also significantly reduces the number of LU
decompositions required for the contour integral.
This makes our algorithm suitable for finding eigenvalues lying within certain
contours for large, sparse problems.

The remainder of this paper is organized as follows.
In Section~\ref{sec:pre}, we provide a review of Beyn's algorithm and the infinite GMRES,
along with some of their useful properties.
In Section~\ref{sec:alg} we introduce the weighting technique and the
two-level orthogonalization technique, and propose the structure of our
algorithm.
Implementation details, including the selection of some parameters and
corresponding analysis, will be discussed in Section~\ref{sec:det}.
Finally, numerical experiments will be presented in Section~\ref{sec:numerexp}
to illustrate the efficiency of our algorithm.

\section{Preliminaries}
\label{sec:pre}

\subsection{Nonlinear eigenvalue problems}
In this work, we focus on nonlinear eigenvalue problems~\eqref{eq:nep}, where
\(T\colon\Omega\rightarrow\mathbb{C}^{n\times n}\) is holomorphic in the
domain \(\Omega\subset\mathbb C\) with sufficiently smooth boundary.
Our goal is to find all eigenvalues \(\lambda_1\), \(\dotsc\), \(\lambda_k\)
lying in \(\Omega\), as well as their corresponding (right) eigenvectors
\(v_1\), \(\dotsc\), \(v_k\).

We say \(\lambda_j\) is a \emph{simple} eigenvalue if
\(\ker\bigl(T(\lambda_j)\bigr)=\Span\{v_j\}\) while
\(T^{\prime}(\lambda_j)v_j\notin\Range\bigl(T(\lambda_j)\bigr)\).
For convenience, in the following, we always assume that there are finitely
many eigenvalues lying in \(\Omega\), and all of them are simple.

\subsection{Beyn's algorithm}
Beyn's algorithm~\cite{Beyn2012} can be regarded as a special case of the
Sakurai--Sugiura algorithm where only the zeroth moment
\[
\mathcal{M}_0=\frac{1}{2\pi\mi}\int_{\partial\Omega}T(\xi)^{-1}Z\md\xi
\]
and first moment
\[
\mathcal{M}_1=\frac{1}{2\pi\mi}\int_{\partial\Omega}\xi T(\xi)^{-1}Z\md\xi
\]
are involved.
For almost any \(Z\in\mathbb{C}^{n\times k}\), if we perform the singular
value decomposition (SVD), \(\mathcal{M}_0=V_0\Sigma_0W_0\herm\),
and construct \(\breve{\mathcal{M}}_1=V_0\herm\mathcal{M}_1W_0\Sigma_0^{-1}\),
it can be proved that \(\lambda_1\), \(\dotsc\), \(\lambda_k\) are exactly
the eigenvalues of \(\breve{\mathcal{M}}_1\).
Furthermore, we can diagonalize \(\breve{\mathcal{M}}_1\) as
\(\breve{\mathcal{M}}_1=S\Lambda S^{-1}\), where
\(\Lambda=\diag\{\lambda_1, \dotsc, \lambda_k\}\) so that
\(V_0S=[v_1, \dotsc, v_k]\) consists of the corresponding eigenvectors
of \(T\).

In Beyn's algorithm, numerical quadrature rules are used to approximate these
moments.
We assume the boundary of the region, \(\partial\Omega\), can be parameterized
as
\[
\varphi\in C^1[0,2\pi],\qquad\varphi(\theta+2\pi)=\varphi(\theta).
\]
Taking \(N\) equidistant quadrature nodes as \(\theta_j=2j\pi/N\) (for \(j=0\),
\(\dotsc\), \(N-1\)) and applying the trapezoidal rule, we will obtain
\begin{equation}
\label{eq:moments}
\begin{aligned}
\mathcal{M}_{0,N}&=\frac{1}{\mi N}\sum_{j=0}^{N-1}\varphi'(\theta_j)
T\bigl(\varphi(\theta_j)\bigr)^{-1}Z,\\
\mathcal{M}_{1,N}&=\frac{1}{\mi N}\sum_{j=0}^{N-1}\varphi(\theta_j)
\varphi'(\theta_j) T\bigl(\varphi(\theta_j)\bigr)^{-1}Z.
\end{aligned}
\end{equation}
In this work, we always use an ellipse contour
\[
\varphi(\theta)=c+a\cos(\theta)+\mi b\sin(\theta)
\]
discretized by the trapezoidal rule with \(N\)
equidistant quadrature nodes.
This is a usual choice in many works~\cite{AST2009,GMP2018,YS2013}, and has
been proved to converge exponentially; see~\cite{TW2014}.

\subsection{Infinite GMRES}
\label{sub-sec:infGMRES}
To approximate the moments by~\eqref{eq:moments}, we have to evaluate
\(T(\xi_j)^{-1}Z\) for several quadrature nodes of the form
\(\xi_j=\varphi(\theta_j)\).
With infinite GMRES (infGMRES)~\cite{CJS2023}, we can solve them efficiently.

Suppose we need to compute \(T(\xi)^{-1}z\) for several values of
\(\xi\)'s around the origin%
\footnote{The origin is chosen for simplicity.
By a simple change of variable \(\xi\gets\xi-\eta\), it is easy to
perform the computation at arbitrary location on the complex plane.}
all at once with infGMRES.
We will firstly approximate \(T\) by a Taylor expansion as
\[
T(\xi)\approx\sum_{j=0}^{p}\frac{\xi^j}{j!}T^{(j)}(0),
\]
and linearize it to \(\mathcal{L}_0-\xi\mathcal{L}_1\), where
\begin{equation}
\label{eq:complin}
\mathcal{L}_0=\bmat{
T(0)&\frac{T^{(1)}(0)}{1!}&\frac{T^{(2)}(0)}{2!}&\cdots&\frac{T^{(p)}(0)}{p!}\\
&I&&&\\
&&I&&\\
&&&\ddots&\\
&&&&I\\
},\quad
\mathcal{L}_1=\bmat{
0&&&&\\
I&0&&&\\
&I&\ddots&&\\
&&\ddots&0&\\
&&&I&0
}.
\end{equation}
It can be proved that the first \(n\) elements of
\((\mathcal{L}_0-\xi\mathcal{L}_1)^{-1}\VECz(z)\) are exactly equal to
\begin{equation}
\label{eq:polysol}
\biggl(\sum_{j=0}^{p}\frac{\xi^j}{j!}T^{(j)}(0)\biggr)^{-1}z,
\end{equation}
where \(\VECz(z)=[z\herm,0,\dotsc,0]\herm\).
Under the assumption that the Taylor expansion is sufficiently accurate,
in order to compute \(T(\xi)^{-1}z\), we just compute
\((\mathcal{L}_0-\xi\mathcal{L}_1)^{-1}\VECz(z)\).

After extracting \(\mathcal{L}_0^{-1}\), a multi-shift GMRES can be employed
to solve for
\((I-\xi\mathcal{L}_1\mathcal{L}_0^{-1})^{-1}\VECz(z)\) for several
\(\xi\)'s easily.
This is because once we obtain
\[
\mathcal{L}_1\mathcal{L}_0^{-1}\mathcal{U}_m=\mathcal{U}_{m+1}\underline{H}_m
\]
by the Arnoldi process, we also have
\[
(I-\xi\mathcal{L}_1\mathcal{L}_0^{-1})\mathcal{U}_m
=\mathcal{U}_{m+1}(\underline{I}_m-\xi\underline{H}_m),
\]
which is exactly the Arnoldi decomposition of
\(I-\xi\mathcal{L}_1\mathcal{L}_0^{-1}\).
Here, \(\mathcal{U}_m\) is the matrix whose columns are Arnoldi vectors,
\(\underline{I}_m\) is the \(m\times m\) identity with an extra zero row at
the bottom, and \(\underline{H}_m\in\mathbb{C}^{(m+1)\times m}\) is an upper
Hessenberg matrix.
Thus, as we already obtained the Hessenberg matrix \(\underline{H}_m\)
by an Arnoldi process on \(\mathcal{L}_1\mathcal{L}_0^{-1}\), we need only to
tackle a small-scale least squares problem
\begin{equation}
\label{eq:ls}
y_{*}=\min_y\bigl\lVert(\underline{I}_m-\xi_{*}\underline{H}_m)y
-\lVert z\rVert_2e_1\bigr\rVert_2,
\end{equation}
for any certain \(\xi=\xi_{*}\).
Since \(\mathcal{U}_my_{*}\) is the approximate value of
\((I-\xi_{*}\mathcal{L}_1\mathcal{L}_0^{-1})^{-1}\VECz(z)\), we take
\(\mathcal{L}_0^{-1}\mathcal{U}_my_{*}\)
as the approximate value of \((\mathcal{L}_0-\xi_{*}\mathcal{L}_1)^{-1}\VECz(z)\).
Then, the first \(n\) elements of \(\mathcal{L}_0^{-1}\mathcal{U}_my_{*}\)
become the approximate value for \(T^{-1}(\xi_{*})z\).
For reference, we list our infGMRES briefly in Algorithm~\ref{alg:infgmres}.
It should be noted that, in the Arnoldi process of Algorithm~\ref{alg:infgmres},
the matrix--vector multiplications with \(\mathcal{L}_0^{-1}\) and
\(\mathcal{L}_1\) can be performed implicitly with matrix--vector
multiplications on \(T^{(j)}(0)\)'s and solving linear systems with \(T(0)\).
The matrices \(\mathcal{L}_0\) and \(\mathcal{L}_1\) are never constructed explicitly.

\begin{algorithm}[tb!]
\caption{infGMRES}
\begin{algorithmic}[1]
\label{alg:infgmres}
\REQUIRE Maximum iteration \(m\), the parameter-dependent matrix
\(T(\xi)\colon\mathbb{C}\rightarrow\mathbb{C}^{n\times n}\), the right-hand
side \(z\in\mathbb{C}^n\) and the points to be solved \(\xi_j\) for \(j=0\),
\(\dotsc\), \(N-1\)
\ENSURE Approximations \(x_{0,j}\approx T(\xi_j)^{-1}z\) for \(j=0\),
\(\dotsc\), \(N-1\)
\STATE Linearize \(T\) to
\[
\mathcal{L}_0=\bmat{
T(0)&\frac{T^{(1)}(0)}{1!}&\cdots&\frac{T^{(p)}(0)}{p!}\\
&I&&\\
&&\ddots&\\
&&&I\\
},\quad
\mathcal{L}_1=\bmat{
0&&&\\
I&0&&\\
&\ddots&\ddots&\\
&&I&0
},\quad p>m
\]
\STATE Perform Arnoldi process on \(\bigl(\mathcal{L}_1\mathcal{L}_0^{-1}, \VECz(z)\bigr)\) to obtain
\(\mathcal{L}_1\mathcal{L}_0^{-1}\mathcal{U}_m=\mathcal{U}_{m+1}\underline{H}_m\)
\STATE Set \(y_j\gets\argmin_y\bigl\lVert(\underline{I}_m-\xi_j\underline{H}_m)y-\lVert z\rVert_2e_1\bigr\rVert_2\) for \(j=0\), \(\dotsc\), \(N-1\)
\STATE Set \(x_{0,j}\gets T(0)^{-1}\bmat{I&-T^{(1)}(0)&\cdots&-T^{(p)}(0)/p!}\mathcal{U}_my_j\) for \(j=0\), \(\dotsc\), \(N-1\)
\label{alg-step:Qnolog}
\end{algorithmic}
\end{algorithm}

It can be seen that
applying \(\mathcal{L}_0^{-1}\) to several vectors involves essentially only
one matrix factorization of \(T(0)\).
Hence, infGMRES is efficient when solving with many \(\xi\)'s.
Furthermore, in~\cite{CJS2023}, it is proved that under certain mild
assumptions the action \(\mathcal{L}_0^{-1}\) can be applied approximately
without affecting the convergence of the algorithm.
This feature makes infGMRES even more attractive, if implemented properly.

Additionally, we note here that the order of the Taylor expansion,
\(p\), does not need to be determined in advance because at the \(j\)th
iteration of GMRES only \(T^{(s)}(0)\)'s with \(s<j\) are involved.
If we process up to some certain iterations, say \(j\), and find that the
desired accuracy is not achieved, we can just provide the algorithm with a new
\(T^{(j+1)}(0)\) and continue.
Therefore, we can always assume that \(p=\infty\), which means the nonlinear
function \(T\) is infinitely approximated by Taylor series
within the domain of convergence.
That is also why these algorithms are usually called infinite, or hold dynamic
polynomial approximation properties~\cite{VMM2013}.
For simplicity, we illustrate the algorithms with the a priori prepared,
order \(p>m\), Taylor expansions in this paper.
But readers should keep in mind that the order \(p\) can be increased during
the algorithm process and does not need to be decided in advance.

\section{Proposed algorithm}
\label{sec:alg}

With the methods we mentioned above, our idea is to solve linear systems in
Beyn's method by infGMRES.
A framework is summarized in Algorithm~\ref{alg:beyn}.
In the following paragraphs, we will discuss how to employ infGMRES
economically and efficiently in Step~\ref{alg-step:infgmres} of
Algorithm~\ref{alg:beyn}.

\begin{algorithm}[tb!]
\caption{Beyn's method with infGMRES}
\begin{algorithmic}[1]
\label{alg:beyn}
\REQUIRE The parameter-dependent matrix
\(T(\xi)\colon\mathbb{C}\rightarrow\mathbb{C}^{n\times n}\), the initial
guess \(Z=[z_1,\dotsc,z_k]\in\mathbb{C}^{n\times k}\), the contour
\(\varphi\) and quadrature nodes \(\theta_j\) for \(j=0\), \(\dotsc\), \(N-1\)
\ENSURE Approximate eigenvalues \(\Lambda\) and eigenvectors \(V\)
\FOR{\(s=1\), \(\dotsc\), \(k\)}
  \STATE Use infGMRES to solve for \(T\bigl(\varphi(\theta_j)\bigr)^{-1}z_s\) for \(j=0\), \(\dotsc\), \(N-1\) simultaneously
  \label{alg-step:infgmres}
\ENDFOR
\STATE Set \(\mathcal{M}_{0,N}\gets\frac{1}{\mi N}\sum_{j=0}^{N-1}\varphi'(\theta_j)T\bigl(\varphi(\theta_j)\bigr)^{-1}Z\)
\STATE Set \(\mathcal{M}_{1,N}\gets\frac{1}{\mi N}\sum_{j=0}^{N-1}\varphi(\theta_j)\varphi'(\theta_j)T\bigl(\varphi(\theta_j)\bigr)^{-1}Z\)
\STATE Singular value decomposition \(\mathcal{M}_{0,N}=V_0\Sigma_0 W_0\herm\)
\STATE Set \(\breve{\mathcal{M}}_{1,N}\gets V_0\herm\mathcal{M}_{1,N}W_0\Sigma_0^{-1}\)
\STATE Eigenvalue decomposition \(\breve{\mathcal{M}}_{1,N}=S\Lambda S^{-1}\)
\STATE Set \(V\gets V_0S\)
\end{algorithmic}
\end{algorithm}

\begin{figure}[tb!]
\centering
\includegraphics[height=6.2cm]{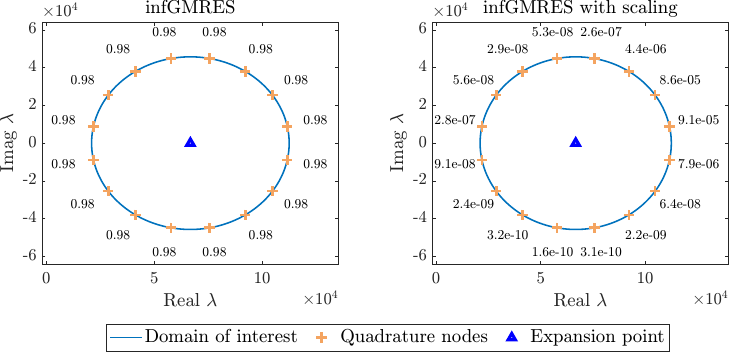} 
\caption{
We use infGMRES to solve several linear systems of the {\tt{gun}} problem.
The contour is a circle centered at \(66762\) with a radius \(45738\).
The quadrature nodes \(\xi_j\) from \(16\) linear systems, \(T(\xi_j)x=z\),
\(j=0\), \(\dotsc\), \(15\), lie equidistantly on the contour.
The only Taylor expansion point is located at the center of
the contour, and the infGMRES method with \(m=32\) maximum iterations is used
in both figures.
To illustrate the accuracy, we plot the relative residuals
\(\lVert T(\xi_j)x_{0,j}-z\rVert_2/\bigl(\lVert T(\xi_j)\rVert_2\lVert x_{0,j}\rVert_2
+\lVert z\rVert_2\bigr)\), where \(x_{0,j}\) stands for the approximate solution of the
\(j\)th linear system (same terms as in Algorithm~\ref{alg:infgmres}).
In practice, for approximating eigenpairs, we usually need these linear
systems to be solved to an accuracy higher than \(10^{-10}\).
The original infGMRES method (left) failed at all \(16\) points,
whereas the infGMRES applied to the variable-substituted system
\(T(5a\tilde\xi+c)\)~(right), although still not accurate enough,
demonstrates a seemingly improved performance.}
\label{fig:gunbad}
\end{figure}

One possible approach is to expand
the Taylor series of \(T\) on the center of the
contour to approximate all the linear systems on the contour at once.
Unfortunately, the original implementation of infGMRES is usually not
accurate enough, especially when the ellipse is relatively large
or some eigenvalues lie close to the quadrature nodes~\cite{CJS2023}.
To illustrate this, we take the {\tt{gun}} problem from the NLEVP
collection~\cite{BHM2013} as an example.
We apply infGMRES on the center of a circular contour and solve the linear
system at each quadrature node; see
Figure~\ref{fig:gunbad}~(left).
With \(32\) GMRES iterations, we obtained completely incorrect solutions. 

Though this failure is partly caused by using too few iterations,
it is not unavoidable within the same number of iterations.
If we instead run infGMRES on the variable-substituted matrix
\(T(5a\tilde\xi+c)\) to solve the same linear systems,
the accuracy can be improved apparently within the same number of iterations;
see Figure~\ref{fig:gunbad}~(right).
Here, we take \(a=45738\) as the radius of the circular contour
and \(c=66762\) as the center of the contour.
We first provide a more general framework of this scaling technique in
Section~\ref{sec:weight}.
Detailed analysis will be included in Section~\ref{sec:weightdetail}.

We remark that, even with the scaling technique, it is impractical to
solve the linear systems for all quadrature nodes \(\xi_j\) on the contour with
a single Taylor expansion at the center of the contour \(c\).
Actually, it is proved in~\cite{JC2022} that the accuracy of the infGMRES
decays with the distance to the Taylor expansion point
\(\lvert \xi_{*}-c\rvert\).
Therefore, when employing infGMRES in contour integral-based algorithms, we
usually expand the Taylor series at several different points to guarantee
that every linear system on the contour can be solved accurately.
These points, which are referred to as \textit{the expansion points}, will be
discussed in Section~\ref{subsec:selectexp}.

The other problem is the redundant memory usage of \(\mathcal{U}_m\).
Storing the whole \(\mathcal{U}_m\), as classical GMRES does,
requires \(\mathcal{O}(m^2n)\) memory.
This will greatly limit the maximum iteration count GMRES can take,
especially when \(n\) is very large.
Following~\cite{KR2014}, we will apply a TOAR-like technique~\cite{LSB2016}
to compress the memory footprint to \(\mathcal{O}(mn+m^3)\).
We will discuss this technique in more detail in Section~\ref{subsec:toar}.

\subsection{Weighting}
\label{sec:weight}
In~\cite[Remark 6.2]{CJS2023}, it has been mentioned that an appropriate
scaling can accelerate the convergence of infGMRES.
But no proof or theoretical analysis is provided to justify the essential
cause.
In this subsection, we provide a more general implementation of the
scaling.
Further insights on this technique and its specific application for
accelerating the convergence of infGMRES will be covered in
Section~\ref{sec:weightdetail}.

When using infGMRES to solve for \(T(\xi_{*})^{-1}z\) with a
Taylor expansion centered at \(0\), we are actually
solving for
\[
(\mathcal{L}_0-\xi_{*}\mathcal{L}_1)^{-1}\VECz(z)
=\bmat{
T_0&T_1&\cdots&T_p\\
-\xi_{*} I&I&&\\
&\ddots&\ddots&\\
&&-\xi_{*} I&I
}^{-1}
\bmat{
z\\
0\\
\vdots\\
0
},
\]
where for simplicity we denote \(T_j=T^{(j)}(0)/j!\), \(j=0\), \(\dotsc\), \(N-1\).
To apply a scaling here, firstly notice that solving for
\(T(\xi)^{-1}z\) at \(\xi=\xi_{*}\) is equivalent to solving for
\(\tilde T(\tilde\xi)^{-1}z\) at \(\tilde\xi=\xi_{*}/\rho\), where
\(\tilde T(\tilde\xi)=T(\rho\tilde\xi)\).
Repeating the same linearization process on \(\tilde T(\tilde\xi)\) yields
alternative companion matrices of the form
\begin{equation}
\label{eq:complinsca}
\tilde{\mathcal{L}}_0=\bmat{
T_0&\rho T_1&\rho^2T_2&\cdots&\rho^pT_p\\
&I&&&\\
&&I&&\\
&&&\ddots&\\
&&&&I
},\quad\tilde{\mathcal{L}}_1=\mathcal{L}_1.
\end{equation}
To solve the original system at \(\tilde\xi=\xi_{*}/\rho\), we solve for
\[
\Bigl(\tilde{\mathcal{L}}_0-\frac{\xi_{*}}{\rho}\tilde{\mathcal{L}}_1\Bigr)^{-1}\VECz(z)
=\bmat{
T_0&\rho T_1&\cdots&\rho^pT_p\\
-\frac{\xi_{*}}{\rho}I&I&&\\
&\ddots&\ddots&\\
&&-\frac{\xi_{*}}{\rho}I&I
}^{-1}
\bmat{
z\\
0\\
\vdots\\
0
}.
\]

Comparing~\eqref{eq:complin} and~\eqref{eq:complinsca}, we notice that
\[
\tilde{\mathcal{L}}_0=D_\rho^{-1}\mathcal{L}_0D_\rho,\quad
\frac{1}{\rho}\tilde{\mathcal{L}}_1=D_\rho^{-1}\mathcal{L}_1D_\rho,\quad
\tilde{\mathcal{L}}_0-\frac{\xi_{*}}{\rho}\tilde{\mathcal{L}}_1=
D_\rho^{-1}(\mathcal{L}_0-\xi_{*}\mathcal{L}_1)D_\rho,
\]
where
\[
D_\rho=
\bmat{
I&&&\\
&\rho I&&\\
&&\ddots&\\
&&&\rho^pI\\
}.
\]
This reminds us the balancing techniques~\cite{CD2000} for eigenvalue
problems, where we may take
\begin{equation}
\label{eq:balancematrix}
D=\bmat{
d_0I&&&\\
&d_1I&&\\
&&\ddots&\\
&&&d_pI\\
}
\end{equation}
with arbitrary \(d_j\)'s.
In fact, we can also obtain the approximate solution~\eqref{eq:polysol}
by solving for
\(\bigl(D^{-1}(\mathcal{L}_0-\xi\mathcal{L}_1)D\bigr)^{-1}\VECz(z)\);
see Lemma~\ref{lem:eqdcomp}.

\begin{lemma}
\label{lem:eqdcomp}
Suppose \(d_j\in\mathbb{C}\backslash\{0\}\) and
\(T_j\in\mathbb{C}^{n\times n}\) for \(j=0\), \(\dotsc\), \(p\), and
\[
D=\bmat{
d_0I&&&\\
&d_1I&&\\
&&\ddots&\\
&&&d_pI\\
},
\]
\[
\mathcal{L}_0=\bmat{
T_0&T_1&T_2&\cdots&T_p\\
&I&&&\\
&&I&&\\
&&&\ddots&\\
&&&&I\\
},\quad
\mathcal{L}_1=\bmat{
0&&&&\\
I&0&&&\\
&I&\ddots&&\\
&&\ddots&0&\\
&&&I&0
}.
\]
Then, for any scalar \(\xi\in\mathbb{C}\) and vector \(z\in\mathbb{C}^n\),
\[
\bigl(D^{-1}(\mathcal{L}_0-\xi\mathcal{L}_1)D\bigr)^{-1}\VECz(z)
=\bmat{
\bigl(\sum_{j=0}^{p}\xi^jT_j\bigr)^{-1}z\\
\ast\\
\vdots\\
\ast
}.
\]
\end{lemma}
\begin{proof}
The proof of Lemma~\ref{lem:eqdcomp} follows from~\cite[Theorem 3.1]{CJS2023}.
We do not repeat it here.
\end{proof}
Lemma~\ref{lem:eqdcomp} indicates that we can perform infGMRES on any balanced
companion linearization \((D^{-1}\mathcal{L}_0D,D^{-1}\mathcal{L}_1D)\),
and a scaling can be regarded as a special case where \(D=D_\rho\).
Additionally, since
\(D^{-1}\mathcal{L}_jD=(d_0^{-1}D)^{-1}\mathcal{L}_j(d_0^{-1}D)\),
\(j=0\), \(1\),
we can always assume \(d_0=1\) without loss of generality.

\subsection{Compact representation of the Arnoldi basis}
\label{subsec:toar}
The two-level orthogonal Arnoldi (TOAR) procedure, initially outlined
in~\cite{ZS2013} (where it is referred to as a compact Arnoldi decomposition),
is a memory-efficient algorithm for computing an orthogonal basis of a
second-order Krylov subspace, which is widely used for reducing the memory
footprint when solving quadratic eigenvalue problems (QEP).
The most remarkable insight of TOAR is that the Arnoldi vectors of a
companion linearization of the QEP can be represented in terms of a common
basis for both upper and lower halves.
The TOAR procedure is proved to be numerically stable~\cite{LSB2016}, and a
framework is developed to extend this technique to more general
cases~\cite{VMM2015}.
In~\cite{JMR2017}, a similar TIAR procedure is implemented on the infinite Arnoldi
algorithm from the point of view of a tensor representation.
In this work, we extend the idea of two-level orthogonalization to reduce the increasingly
larger memory usage.

Firstly, partition the Arnoldi vectors \(\mathcal{U}_m\) by row blocks to
\[
\mathcal{U}_m=\bmat{
U_{m,0}\\
\vdots\\
U_{m,p}},\quad U_{m,j}\in\mathbb{C}^{n\times m},\quad j=0,\dotsc,p.
\]
Then, there is a \(Q_m\in\mathbb{C}^{n\times m}\), \(Q_m\herm Q_m=I\),
such that \(U_{m,j}=Q_m\breve U_{m,j}\) for some
\(\breve U_{m,j}\in\mathbb{C}^{m\times m}\) (for \(j=0\), \(\dotsc\), \(p\)).
Furthermore, the Arnoldi process can be rearranged into
\[
\mathcal{L}_1\mathcal{L}_0^{-1}\mathcal{U}_m
=\mathcal{L}_1\mathcal{L}_0^{-1}
\bmat{
Q_m\breve U_{m,0}\\
\vdots\\
Q_m\breve U_{m,p}}
=\bmat{
Q_m\breve U_{m,0}&[Q_m,q_{m+1}]\breve u_{m+1,0}\\
\vdots&\vdots\\
Q_m\breve U_{m,p}&[Q_m,q_{m+1}]\breve u_{m+1,p}}
\underline{H}_m
=\mathcal{U}_{m+1}\underline{H}_m,
\]
where \(\underline{H}_m\in\mathbb{C}^{(m+1)\times m}\), \(\breve u_{m+1,j}\in\mathbb{C}^{m+1}\) (for \(j=0\), \(\dotsc\), \(p\)), and the
columns of both \([Q_m,q_{m+1}]\) and
\[
\bmat{
Q_m\breve U_{m,0}&[Q_m,q_{m+1}]\breve u_{m+1,0}\\
\vdots&\vdots\\
Q_m\breve U_{m,p}&[Q_m,q_{m+1}]\breve u_{m+1,p}}
\]
are orthonormal.
Thus, by setting \(Q_{m+1}=[Q_m,q_{m+1}]\) and
\[
\breve U_{m+1,j}
=\bmat{
\begin{tabular}{c|} \(\breve U_{m,j}\) \\ 0 \end{tabular}
& \breve u_{m+1,j-1}
},
\]
the Arnoldi vectors \(\mathcal{U}_{m+1}\) share the same structure
\[
\mathcal{U}_{m+1}
=\bmat{Q_{m+1}\breve U_{m+1,0} \\ \vdots \\ Q_{m+1}\breve U_{m+1,p}}.
\]
Therefore, we only need to store \(Q_m\) and \(\breve U_{m,j}\) in memory for
the Arnoldi process.
Additionally, thanks to the pattern of \(\VECz(z)\), we have
\(\breve U_{m,j}=0\) for \(j\ge m\).
Overall, the memory usage will be \(\mathcal{O}(mn+m^3)\).
However, since \(\mathcal{U}_m\) is not explicitly formed, additional care
should be taken when orthogonalizing, using the so called two-level
orthogonalization.
Combined with weighting, we list our weighted two-level orthogonal
infGMRES in Algorithm~\ref{alg:wtinfgmres} for reference.

\begin{algorithm}
\caption{Weighted two-level orthogonal infGMRES}
\begin{algorithmic}[1]
\label{alg:wtinfgmres}
\REQUIRE Maximum iterations \(m\), the matrices \(T_j\in\mathbb{C}^{n\times n}\) for \(j=0\), \(\dotsc\), \(p\), \(p>m\), the initial guess \(z\in\mathbb{C}^n\) and the weights \(d_1\), \(\dotsc\), \(d_p\) 
\ENSURE A function \(f(\cdot)\) such that \(f(\xi)\) approximates \(T(\xi) ^{-1}z\)
\STATE \textbf{function} \([f]={\tt{WTinfGMRES}}(T_0,\dotsc,T_p,d_1,\dotsc,d_p,z)\)
\STATE \(Q_1\gets z/\lVert z\rVert_2\), \(H_0\gets [~]\), \(Q_{\save}\gets [~]\)
\STATE \(\breve U_{1,0}\gets 1\), \(\breve u_{1,0}\gets 1\), \(\breve U_{1,s}\gets 0\), \(\breve u_{1,s}\gets 0\), \(s=1\), \(\dotsc\), \(p\)
\FOR{\(j=1\), \(\dotsc\), \(m\)}
  \STATE \(q_{j+1}\gets d_1^{-1}T_0^{-1}\bigl(Q_j\breve u_{j,0}-\sum_{s=1}^pd_sT_sQ_j\breve u_{j,s}\bigr)\)
  \label{alg-step:solvem}
  \STATE \(Q_{\save}\gets [Q_{\save}, d_1q_{j+1}]\) \quad\COMMENT{Record \(\mathcal{L}_0^{-1}Q_j\)}
  \label{alg-step:Qlogsave}
  \STATE \(l_j\gets Q_j\herm q_{j+1}\) \quad\COMMENT{First level orthogonalization}
  \STATE \(q_{j+1}\gets q_{j+1}-Q_jl_j\)
  \STATE \(\alpha_j\gets \lVert q_{j+1}\rVert_2\)
  \STATE \(q_{j+1}\gets q_{j+1}/\alpha_j\)
  \STATE \(Q_{j+1}\gets[Q_j,q_{j+1}]\)
  \STATE \(h_j=\breve U_{j,1}\herm l_j+\sum_{s=2}^p(d_{s-1}/d_s)\breve U_{j,s}\herm \breve u_{j,s-1}\) \quad\COMMENT{Second level orthogonalization}
  \STATE \(\breve u_{j+1,0}\gets 0\), \(\breve u_{j+1,1}\gets\bmat{l_j-\breve U_{j,1}h_j\\\alpha_j}\),
  \(\breve u_{j+1,s}\gets\bmat{\breve u_{j,s-1}-\breve U_{j,s}h_j\\0}\), \(s=2\), \(\dotsc\), \(p\)
  \STATE \(\beta_j\gets\sqrt{\sum_{s=0}^j\breve u_{j+1,s}\herm \breve u_{j+1,s}}\)
  \STATE \(\breve u_{j+1,s}\gets \breve u_{j+1,s}/\beta_j\), \(s=0\), \(\dotsc\), \(p\)
  \STATE \(H_j\gets\bmat{H_{j-1}&h_j\\0&\beta_j}\) \quad\COMMENT{Update \(U\) and \(H\)}
  \STATE \(\breve U_{j+1,s}\gets\bmat{\begin{tabular}{c}\(\breve U_{j,s}\)\\\(0\)
  \end{tabular}&\breve u_{j+1,s}}\), \(s=0\), \(\dotsc\), \(p\)
\ENDFOR
\STATE \textbf{return} \(f\colon \xi\mapsto\lVert z\rVert_2Q_{\save}(I-\xi H_m)^{\dagger}e_1\)
\label{alg-step:Qloguse}
\STATE \textbf{end function}
\end{algorithmic}
\end{algorithm}

\begin{remark}
For simplicity and consistency, we initialize \(\breve U_{j,s}\)'s for
\(s=0\), \(\dotsc\), \(p\) from the beginning.
However, readers should remember that \(\breve U_{j,s}=0\) for \(s\ge j\).
Hence, \(\breve U_{j,s}\)'s are neither computed nor stored in memory.
\end{remark}

\begin{remark}
Notice that in Step~\ref{alg-step:Qnolog} of Algorithm~\ref{alg:infgmres},
we have to apply
\[
Q_{\save}=T(0)^{-1}\bmat{I&-T^{(1)}(0)&\cdots&-T^{(p)}(0)/p!}\mathcal{U}_m
\]
on each \(y_j\) to recover the approximate solution \(x_{0,j}\),
which may involve extra workload.
Fortunately, as mentioned in~\cite{CJS2023},
we can store \(Q_{\save}\) column by column during the infGMRES process,
to make this cheaper; see Step~\ref{alg-step:Qlogsave} and
Step~\ref{alg-step:Qloguse} in Algorithm~\ref{alg:wtinfgmres}.
\end{remark}

\begin{remark}
We remark here that even with the TOAR-like technique,
\(Q_m\) and
\[
\breve{\mathcal{U}}_m=\bmat{\breve U_{m,0}\\
\vdots\\
\breve U_{m,p}
}
\]
may not be fully orthogonal in finite precision arithmetic.
As we show in Table~\ref{tab:ortho}, Algorithm~\ref{alg:wtinfgmres} loses the
orthogonality of \(Q_m\) in almost every test example and sometimes loses the
orthogonality of \(\breve{\mathcal{U}}_m\) as well.
However, the loss of orthogonality does not always destroy the convergence of
GMRES; see~\cite{GRS1997}.
In fact, we can perform a reorthogonalization on \(Q_m\) and
\(\breve{\mathcal{U}}_m\) to make them orthogonal; see Table~\ref{tab:ortho}.
Nevertheless, the accuracy of infGMRES will not increase, even though both
\(Q_m\) and \(\breve{\mathcal{U}}_m\) are numerically orthogonal.
Therefore, we will not use reorthogonalization in our numerical experiments.
\end{remark}

\begin{table}[tb!]
\centering
\caption{Application of Algorithm~\ref{alg:wtinfgmres} with/without
reorthogonalization to all \(9\) test examples from Table~\ref{tab:expsinfo}.
The last two columns show the orthogonality of \(Q_m\) and
\(\breve{\mathcal{U}}_m\).
The test is performed on all the expansion points, and only the worst values
are shown here.}
\begin{NiceTabular}{cccc}
\hline
&Problem & \(\lVert Q_m\herm Q_m-I\rVert_2\) & \(\lVert \breve{\mathcal{U}}_m\herm \breve{\mathcal{U}}_m-I\rVert_2\)\\
\hline
\multirow{9}*{w/o reorthog.}
&{\tt{spring}}             & \(2.3\times 10^{1}\) & \(6.3\times 10^{-7}\) \\
&{\tt{acoustic\_wave\_2d}} & \(2.1\times 10^{-8}\)  & \(5.8\times 10^{-11}\) \\
&{\tt{butterfly}}          & \(7.8\times 10^{-4}\)  & \(1.4\times 10^{-3}\)  \\
&{\tt{loaded\_string}}     & \(1.6\times 10^{1}\)  & \(1.5\times 10^{-13}\)  \\
&{\tt{photonics}}           & \(1.9\times 10^{1}\)  & \(7.1\times 10^{-11}\)  \\
&{\tt{railtrack2\_rep}}    & \(1.7\times 10^{1}\)  & \(4.0\times 10^{-14}\)  \\
&{\tt{hadeler}}            & \(1.2\times 10^{-14}\)  & \(2.0\times 10^{-14}\)   \\
&{\tt{gun}}                & \(1.8\times 10^{1}\)  & \(4.0\times 10^{-14}\)   \\
&{\tt{canyon\_particle}}   & \(1.9\times 10^{1}\)  & \(8.1\times 10^{-6}\)   \\
\hdottedline
\multirow{9}*{with reorthog.}
&{\tt{spring}}             & \(1.1\times 10^{-15}\)  & \(5.0\times 10^{-16}\)  \\
&{\tt{acoustic\_wave\_2d}} & \(6.8\times 10^{-16}\)  & \(5.0\times 10^{-16}\) \\
&{\tt{butterfly}}          & \(6.3\times 10^{-16}\)  & \(6.7\times 10^{-16}\)  \\
&{\tt{loaded\_string}}     & \(1.1\times 10^{-15}\)  & \(7.0\times 10^{-16}\)  \\
&{\tt{photonics}}           & \(7.7\times 10^{-16}\)  & \(8.1\times 10^{-16}\)  \\
&{\tt{railtrack2\_rep}}    & \(5.7\times 10^{-16}\)  & \(4.3\times 10^{-16}\)  \\
&{\tt{hadeler}}            & \(9.7\times 10^{-16}\)  & \(3.3\times 10^{-16}\)   \\
&{\tt{gun}}                & \(9.1\times 10^{-16}\)  & \(7.4\times 10^{-16}\)   \\
&{\tt{canyon\_particle}}   & \(9.6\times 10^{-16}\)  & \(6.7\times 10^{-16}\)   \\
\hline
\end{NiceTabular}
\label{tab:ortho}
\end{table}

\section{Implementation details of the algorithm}
\label{sec:det}

While our modified infGMRES has been fully described in previous sections,
the implementation process involves determining numerous parameters,
including the scalars \(d_j\), the locations of expansion points \(\eta_j\), and
others.
In this section, we will reveal the connections among these parameters and
the convergence of infGMRES, offering practical guidance on their selection,
and hence finishing the last details needed in Algorithm~\ref{alg:beyn}.

\subsection{Scaling is essentially a weighting strategy on GMRES}
\label{sec:weightdetail}
In Section~\ref{sec:weight}, we have already showed that a scaling can be
generalized into a balanced companion linearization, say
\(D^{-1}\mathcal{L}_0D\) and \(D^{-1}\mathcal{L}_1D\).
Here, we declare that this technique is essentially a weighting strategy for
the least squares problems within GMRES, which is very similar to the
weighted GMRES algorithm~\cite{Essai1998}.
If used appropriately, this technique can help to find a more accurate
solution in a certain search space.

To see this, we have to look into the process of GMRES.
First notice that if the Krylov subspace
\(\mathcal{K}_m\bigl(I-\xi_{*}\mathcal{L}_1\mathcal{L}_0^{-1},\VECz(z)\bigr)\)
is spanned by \(\mathcal{U}_m\), by induction, we can prove that
\(\mathcal{K}_m\bigl(D^{-1}(I-\xi_{*}\mathcal{L}_1\mathcal{L}_0^{-1})D,\VECz(z)\bigr)\)
is spanned by \(D^{-1}\mathcal{U}_m\).
Thus, at the \(m\)th iteration, infGMRES will give the solution by
\begin{equation}
\label{eq:gmresy}
\begin{aligned}
y_{*}&=\argmin_y
\left\lVert
D^{-1}(I-\xi_{*}\mathcal{L}_1\mathcal{L}_0^{-1})D
D^{-1}\mathcal{U}_my-\VECz(z)\right\rVert_2\\
&=\argmin_y\left\lVert
D^{-1}\bigl((I-\xi_{*}\mathcal{L}_1\mathcal{L}_0^{-1})
\mathcal{U}_my-\VECz(z)\bigr)
\right\rVert_2,
\end{aligned}
\end{equation}
where the last equality follows because \(d_0=1\).

Equation~\eqref{eq:gmresy} provides us with several insights.
Firstly, irrespective of the selection of the weighting matrix \(D\),
the final approximation of \(T(\xi_{*})^{-1}z\) will be selected from the same
search space.
That is because the approximate solution will be the first \(n\) elements of
\(D^{-1}\mathcal{L}_0^{-1}\mathcal{U}_my_{*}\) and \(d_0=1\).
The role \(D\) plays in infGMRES is actually a weighting matrix in the least
squares problems~\eqref{eq:gmresy}.
We can choose \(D\) appropriately to guide GMRES to pick a more accurate
solution from \(\mathcal{U}_m\).

The previous sentence may seem contradictory because GMRES already provides
the best solution that minimizes the residual norm.
However, what we really care about is not the residual of
GMRES~\eqref{eq:gmresy}, but
\[
\lVert r_{\rm{N}}\rVert_2=\bigl\lVert T(\xi_{*})x_{0,*}-z\bigr\rVert_2,
\]
where \(x_{0,*}\) is the approximate value of
\(T(\xi_{*})^{-1}z\).
Since
\[
\lVert r_{\rm{N}}-r_{\rm{P}}\rVert_2
=\Bigl\lVert\frac{T^{(p+1)}(\zeta)}{(p+1)!}\xi_{*}^{n+1}\Bigr\rVert_2
\le\frac{\lVert T^{(p+1)}(\zeta)\rVert_2}{(p+1)!}\lVert\xi_{*}\rVert_2^{n+1},
\quad\zeta=\delta\xi_{*},
\quad\delta\in(0,1),
\]
we can instead consider
\begin{equation}
\label{eq:polyres}
\lVert r_{\rm{P}}\rVert_2=\Bigl\lVert\sum_{j=0}^p\xi_{*}^{j}T_jx_{0,*}-z\Bigr\rVert_2
\end{equation}
when \(\mathcal{O}\bigl(\lVert T^{(p+1)}(\zeta)\rVert_2
\lVert\xi_{*}\rVert_2^{n+1}/(p+1)!\bigr)\) is negligible.
In an ideal scenario, we choose the weighting matrix \(D\)
so that~\eqref{eq:gmresy} returns the best solution in the sense
of~\eqref{eq:polyres}.
Therefore, it is crucial to reveal the relationship between~\eqref{eq:gmresy}
and~\eqref{eq:polyres}.

\begin{lemma}
\label{lem:relares}
Suppose \(z\), \(\mathcal{U}_m\), \(\mathcal{L}_0\), \(\mathcal{L}_1\)
and \(T_j\) (for \(j=0\), \(\dotsc\), \(p\)) are all defined as before.
The vector \(y_{*}\in\mathbb{C}^m\) is the solution of~\eqref{eq:gmresy}
for some weighting matrix \(D\).
Then, the polynomial-wise residual
\(r_{\rm{P}}=\sum_{j=0}^p\xi_{*}^{j}T_jx_{0,*}-z\) can be represented as
\begin{equation}
\label{eq:resequal}
r_{\rm{P}}=\biggl[I,-\sum_{j=1}^p\xi_{*}^{j-1}T_{j},-\sum_{j=2}^p\xi_{*}^{j-2}T_j,\dotsc,-T_p\biggr]
\bmat{r_0\\
r_1\\
\vdots\\
r_p},
\end{equation}
where
\begin{equation}
\label{eq:orires}
\bmat{r_0\\
r_1\\
\vdots\\
r_p}=
(I-\xi_{*}\mathcal{L}_1\mathcal{L}_0^{-1})
\mathcal{U}_my_{*}-\VECz(z).
\end{equation}
\end{lemma}

\begin{proof}
We know that if \(y_{*}\) is the solution of~\eqref{eq:gmresy}, infGMRES will
give \(x_{*}=D^{-1}\mathcal{L}_0^{-1}\mathcal{U}_my_{*}\) as the approximate
value of 
\(\bigl(D^{-1}(\mathcal{L}_0-\xi_{*}\mathcal{L}_1)D\bigr)^{-1}\VECz(z)\).
Then, by Lemma~\ref{lem:eqdcomp}, the first \(n\) elements of \(x_{*}\),
denoted by \(x_{0,*}\in\mathbb{C}^n\), will be taken as the approximate
value of \(T(\xi_{*})^{-1}z\).
In other words, we have
\[
D^{-1}\mathcal{L}_0^{-1}\mathcal{U}_my_{*}=x_{*}=
\bmat{
x_{0,*}\\
x_{1,*}\\
\vdots\\
x_{p,*}
},
\]
where \(x_{0,*}\) here is same as the one in~\eqref{eq:polyres}.

Therefore, we can represent~\eqref{eq:orires} in terms of \(x_{*}\) instead
of \(y_{*}\):
\begin{equation}
\label{eq:gmresres}
\bmat{r_0\\
r_1\\
\vdots\\
r_p}
=(I-\xi_{*}\mathcal{L}_1\mathcal{L}_0^{-1})\mathcal{U}_my_{*}-\VECz(z)
=(\mathcal{L}_0-\xi_{*}\mathcal{L}_1)Dx_{*}-\VECz(z).
\end{equation}
Just listing out the equations in~\eqref{eq:gmresres}
\[
\begin{aligned}
d_0T_0x_{0,*}+d_1T_1x_{1,*}+\cdots+d_pT_px_{p,*}&=r_0+z,\\
-\xi_{*} d_0x_{0,*}+d_1x_{1,*}&=r_1,\\
&\vdots\\
-\xi_{*} d_{p-1}x_{p-1,*}+d_px_{p,*}&=r_p,\\
\end{aligned}
\]
and substituting all other equations into the first one
yields~\eqref{eq:resequal}.
\end{proof}

Remember that the solution \(y_{*}\) minimizes~\eqref{eq:gmresy},
so that \(r_j\)'s minimize
\[
\left\lVert
D^{-1}\bigl((I-\xi_{*}\mathcal{L}_1\mathcal{L}_0^{-1})
\mathcal{U}_my_*-\VECz(z)\bigr)
\right\rVert_2=\left\lVert
D^{-1}\bmat{r_0\\
r_1\\
\vdots\\
r_p}
\right\rVert_2=
\sqrt{\sum_{j=0}^p\frac{\lVert r_j\rVert_2^2}{d_j^2}},
\]
which means we can adjust the magnitude of \(r_j\) by setting \(d_j\)
properly.
If we set \(d_j\) to be larger,
then, \(\lVert r_j\rVert_2\) must be more modest.

On the other hand, we know from Lemma~\ref{lem:relares} that
\[
\lVert r_{\rm{P}}\rVert_2\le
\lVert r_0\rVert_2+\Bigl\lVert\sum_{j=1}^{p}\xi_{*}^{j-1}T_j\Bigr\rVert_2\lVert r_1\rVert_2+\cdots+
\lVert T_p\rVert_2\lVert r_p\rVert_2.
\]
Therefore, intuitively, if we make \(\lVert r_j\rVert_2\) relatively small
for larger \(\lVert\sum_{j=s}^p\xi_{*}^{j-s}T_j\rVert_2\),
and relatively large for smaller \(\lVert\sum_{j=s}^p\xi_{*}^{j-s}T_j\rVert_2\),
overall, \(\lVert r_{\rm{P}}\rVert_2\) could be modest.

With this insight, we may set the weights as \(d_0=1\) and
\(d_s=\lVert\sum_{j=s}^p\xi_{*}^{j-s}T_j\rVert_2^{-1}\) for \(s=1\), \(\dotsc\), \(p\).
However, there are cases where \(\lVert\sum_{j=s}^p\xi_{*}^{j-s}T_j\rVert_2\ll1\)
or \(\lVert\sum_{j=s}^p\xi_{*}^{j-s}T_j\rVert_2\gg1\),
which will make \(D\) extremely ill-conditioned.
Thus, in practice, we prefer to balance \(d_j\) by
\[
d_j\gets d_j\cdot\frac{d_2}{d_1^2},\quad j>0.
\]
Consequently, for a certain quadrature node \(\xi_{*}\), a reasonable choice of weights \(d_j\) is:
\[
d_0=1,\quad
d_s=\frac{\gamma}{\lVert\sum_{j=s}^p\xi_{*}^{j-s}T_j\rVert_2},
\quad\gamma=\frac{\lVert\sum_{j=1}^p\xi_{*}^{j-2}T_j\rVert_2^2}{\lVert\sum_{j=2}^p\xi_{*}^{j-3}T_j\rVert_2},
\quad s=1,\dotsc,p.
\]

\begin{remark}[How to determine \(d_j\)'s in practice]
It is impractical to assign a different weighting matrix \(D\) to each
quadrature node, since a different \(D\) results in a different Arnoldi
process, which is very expensive and contradicts our original intention.
A wiser strategy is to use the same \(D\) for all the quadrature nodes under
a given expansion point \(\eta_t\).
In practice, we set the weights as
\begin{equation}
\label{eq:ourwt}
\begin{aligned}
&d_0=1,\quad d_s=\frac{\gamma}{\lVert\sum_{j=s}^p\nu^{j-s}T_j\rVert_2},
\quad\gamma=\frac{\lVert\sum_{j=1}^p\nu^{j-2}T_j\rVert_2^2}{\lVert\sum_{j=2}^p\nu^{j-3}T_j\rVert_2},
\quad s=1,\dotsc,p,\\
\end{aligned}
\end{equation}
where
\[
\nu=2\cdot\max_{j\in\Psi}\lvert\xi_j-\eta_t\rvert,
\qquad\Psi=\bigl\{j\colon\lvert\xi_j-\eta_t\rvert<\lvert\xi_j-\eta_s\rvert,
~s\neq t\bigr\}.
\]
This choice performs well in numerical experiments.
\end{remark}

\begin{remark}
Our analysis does not depend on specific styles of the companion
linearizations.
For example, if we use the classical companion linearization
\[
\mathcal{L}_0=\bmat{
A_1&A_2&\cdots&A_p\\
I&&&\\
&\ddots&&\\
&&I&
},\qquad
\mathcal{L}_1=\bmat{
A_0&\\
&-I\\
&&\ddots&\\
&&&-I
},
\]
it can also be proved that
\[
D_\rho^{-1}(\xi\mathcal{L}_1^{-1}\mathcal{L}_0-I)D_\rho=
\frac{\xi}{\rho}\tilde{\mathcal{L}}_1^{-1}\tilde{\mathcal{L}}_0-I,
\]
where \(\tilde{\mathcal{L}}_0\) and \(\tilde{\mathcal{L}}_1\)
are the companion linearization matrices corresponding to the scaled system \(\tilde T\).
Other analyses on \(D\) will follow.
\end{remark}

\subsection{On polynomial eigenvalue problems}

Since infGMRES is originally designed for non-polynomial eigenvalue problems,
additional care must be taken when dealing with a PEP, or in other words when
\[
T(\xi)=T_0+\xi T_1+\cdots+\xi^gT_g.
\]
In these cases, it is not necessary for us to use infGMRES.
On the contrary, just linearizing it to
\begin{equation}
\label{eq:pepl1}
\breve{\mathcal{L}}_0=\bmat{
T_0&T_1&T_2&\cdots&T_g\\
&I&&&\\
&&I&&\\
&&&\ddots&\\
&&&&I\\
},\quad
\breve{\mathcal{L}}_1=\bmat{
0&&&&\\
I&0&&&\\
&I&\ddots&&\\
&&\ddots&0&\\
&&&I&0
},
\end{equation}
and using multi-shift GMRES to solve for
\((\breve{\mathcal{L}}_0-\xi\breve{\mathcal{L}}_1)^{-1}\VECz(z)\) is
sufficient to solve the problem.

However, if infGMRES is chosen, the linearization becomes
\begin{equation}
\label{eq:pepl2}
\mathcal{L}_0=\bmat{
T_0&T_1&\cdots&T_g&0&\cdots\\
&I&&&\\
&&\ddots&&&\\
&&&I&&\\
&&&&I&\\
&&&&&\ddots
},\quad
\mathcal{L}_1=\bmat{
0&&&&&\\
I&0&&&&\\
&I&\ddots&&&\\
&&\ddots&0&&\\
&&&I&0&\\
&&&&\ddots&\ddots
},
\end{equation}
where we use \(0\) to fill the position of \(T_j\), \(j>g\)
because \(T^{(j)}(\xi)=0\), \(j>g\).
It is interesting to note that \((\mathcal{L}_0-\xi\mathcal{L}_1)^{-1}\VECz(z)\)
also gives the true value of \(T(\xi)^{-1}z\) on the first
\(n\) elements.
A question naturally arises --- is there any performance difference
between~\eqref{eq:pepl1} and~\eqref{eq:pepl2}?

From the point of view of computational cost and memory usage,
in the \(j\)th iteration of~\eqref{eq:pepl1},
the cost of orthogonalization will be \(\mathcal{O}(jn+j^2g)\)
and the total memory usage is also \(\mathcal{O}(jn+j^2g)\).
As for~\eqref{eq:pepl2}, they are \(\mathcal{O}(jn+j^3)\).
For general cases \(m>g\), so that \eqref{eq:pepl1} is cheaper.
However, since \(g<m\ll n\), there is not much difference.

Things become interesting if we look from the perspective of convergence rate.
Equation \eqref{eq:pepl1} can actually be regarded as~\eqref{eq:pepl2}
with the weighting we mentioned in Section~\ref{sec:weightdetail}.
Remember that, with our weighting strategy,
the weighting matrix \(D\) will be
\[
D=\bmat{
d_0I&&&&\\
&d_1I&&&&\\
&&\ddots&&\\
&&&d_gI&\\
&&&&\infty I&\\
&&&&&\ddots
}.
\]
Then, if we regard \(\infty\) as a very large number that can be used in
arithmetic computation, we will have the weighted system
\[
D^{-1}\mathcal{L}_0D=
\bmat{
T_0&T_1&\cdots&T_g&0&\cdots\\
&I&&&\\
&&\ddots&&&\\
&&&I&&\\
&&&&I&\\
&&&&&\ddots
},\]
\[
D^{-1}\mathcal{L}_1D=
\bmat{
0&&&&&\\
\frac{d_0}{d_1}I&0&&&&\\
&\ddots&\ddots&&&\\
&&\frac{d_{g-1}}{d_g}I&0&&\\
&&&0&0&\\
&&&&I&0\\
&&&&&\ddots&\ddots
}.
\]
Notice that \(D^{-1}\mathcal{L}_1\mathcal{L}_0^{-1}D\) is a diagonal block
matrix now.
Since the right-hand side \(\VECz(z)\) has non-zero elements only on the
first \(n\) dimensions, solving for
\(\bigl(D^{-1}(I-\xi\mathcal{L}_1\mathcal{L}_0^{-1})D\bigr)^{-1}\VECz(z)\)
and \(\bigl(D_g^{-1}(I-\xi\breve{\mathcal{L}}_1\breve{\mathcal{L}}_0^{-1})D_g\bigr)^{-1}\VECz(z)\) 
with GMRES will result in exactly the same process, where
\(D_g=\diag\{d_0I,\dotsc,d_gI\}\) is the truncated \(D\).

Since both the cost and the convergence of the two methods are similar,
we also use infGMRES in Section~\ref{sec:numerexp} for PEP to keep results uniform.

\subsection{Selecting multiple expansion points}
\label{subsec:selectexp}
\begin{figure}[tb!]
\centering
\includegraphics[height=6.4cm]{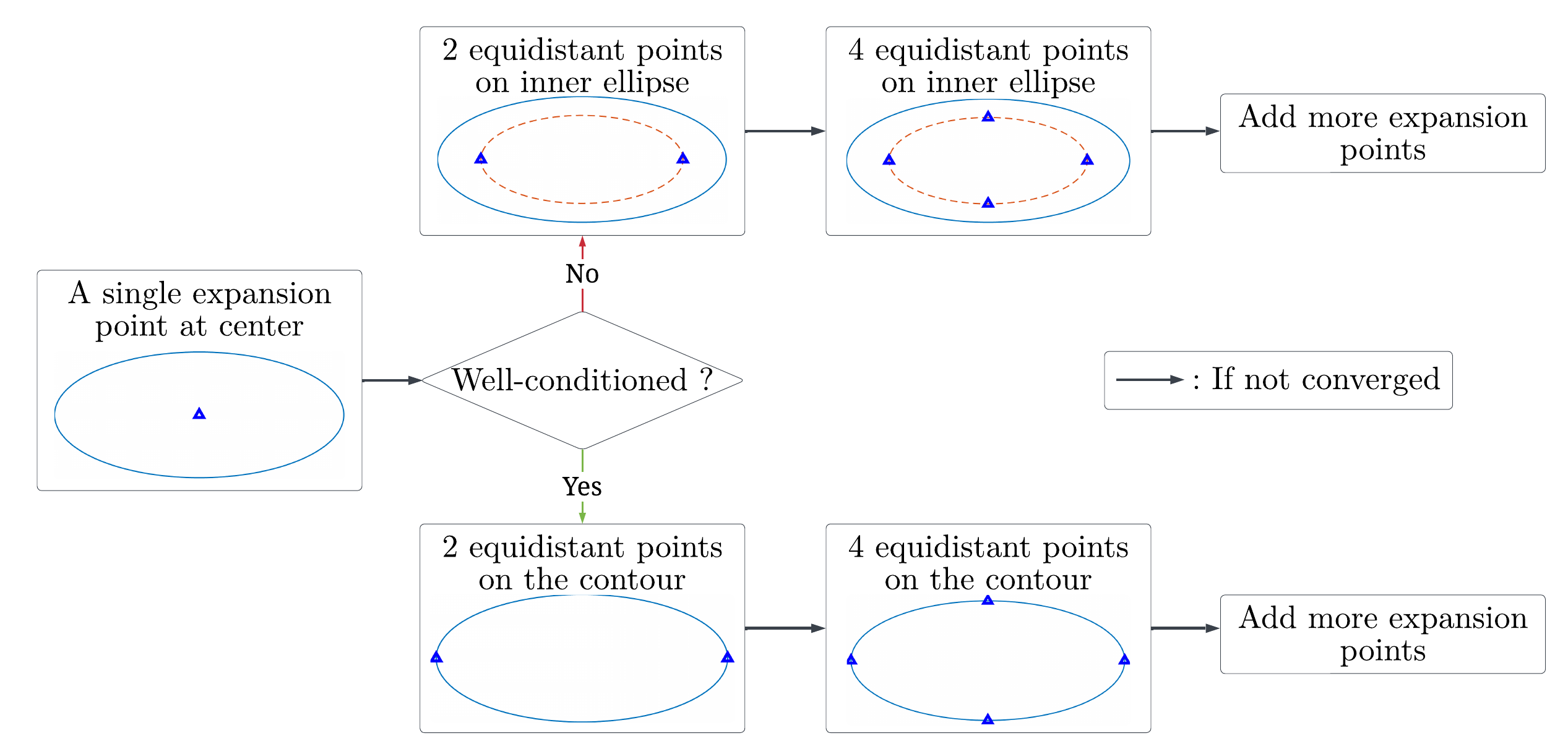}
\caption{A brief guideline on choosing expansion points:
If a single point at the center does not work,
choose expansion points equidistantly on the contour or on an inner ellipse,
based on the condition number of the problem.
Continue doubling the number of expansion points till the accuracy is
satisfactory.}
\label{fig:flowchart}
\end{figure}

Even with the weighting strategy,
it may still take too many iterations for infGMRES to converge,
especially for quadrature nodes that are far from expansion points.
Therefore, it is sensible to use more than one expansion point in practice.
However, a good choice of the expansion points
depends on the singularities, eigenvalues and many other things.
Therefore, adaptively selecting the expansion points is usually impractical,
and it is not the focus of interest in this paper.
In this subsection, we provide a flowchart
for users to heuristically determine the distribution of expansion points;
see Figure~\ref{fig:flowchart}.

At this point, we can fully describe the computation needed for
Step~\ref{alg-step:infgmres} of Algorithm~\ref{alg:beyn};
see Algorithm~\ref{alg:solve}.

\begin{algorithm}
\caption{Step~\ref{alg-step:infgmres} of Algorithm~\ref{alg:beyn}}
\begin{algorithmic}[1]
\label{alg:solve}
\REQUIRE The parameter-dependent matrix
\(T(\xi)\colon\mathbb{C}\rightarrow\mathbb{C}^{n\times n}\),
the right-hand side \(z\in\mathbb{C}^n\), the contour \(\varphi\) and
quadrature nodes \(\theta_j\) for \(j=0\), \(\dotsc\), \(N-1\)
\ENSURE Approximations \(x_{0,j}\approx T\bigl(\varphi(\theta_j)\bigr)^{-1}z\) for \(j=0\), \(\dotsc\), \(N-1\)
\STATE Determine expansion points \(\eta_0\), \(\dotsc\), \(\eta_{\nep}\) as in Figure~\ref{fig:flowchart}
\FOR{\(j=0\), \(\dotsc\), \(\nep\)}
  \STATE Set \(T_{s,j}\gets T^{(s)}(\eta_j)/s!\) for \(s=0\), \(\dotsc\), \(p\)
  \STATE Compute weights \(d_{s,j}\) as in Section~\ref{sec:weight}
  \STATE Set \(f_j(\cdot)\gets{\tt{WTinfGMRES}}(T_{0,j},\dotsc,T_{p,j},d_{0,j},\dotsc,d_{p,j},z)\)~(Algorithm~\ref{alg:wtinfgmres})
\ENDFOR
\FOR{\(j=0\), \(\dotsc\), \(N-1\)}
  \STATE Find \(\eta_s\in\{\eta_0,\dotsc,\eta_{\nep}\}\) closest to \(\varphi(\theta_j)\)
  \STATE Set \(x_{0,j}\gets f_s\bigl(\varphi(\theta_j)\bigr)\)
\ENDFOR
\end{algorithmic}
\end{algorithm}

\begin{remark}[Accuracy of linear systems and the quadrature rule]
In contour integral-based algorithms, the number of quadrature nodes is
sometimes increased in order to improve the accuracy of the quadrature rule.
The question thus arises as to whether the accuracy of the solution of linear
systems should also be increased.
Furthermore, in order to attain a specific level of accuracy for the
quadrature rule, how accurate should the linear systems be solved?
To the best of our knowledge, this is still an open problem.
Nevertheless, in this work, we recommend the users always solve the linear
systems to machine precision.
This is because of the fact that adding quadrature nodes introduces little
additional overhead in our algorithm.
Therefore, a well-distributed expansion point set can be reused whenever the
users need a higher level of accuracy.
\end{remark}
\begin{remark}[Balancing between iterations of infGMRES and number of
expansion points]
Having already addressed the question of the desired level of accuracy for
linear systems, it is now necessary to consider how this accuracy can be
achieved.
Two principal parameters can be adjusted --- the iterations of infGMRES and
the number of expansion points.
In practice, it can be observed that increasing the iterations of infGMRES
sometimes only leads to slight improvements on
the accuracy of solutions of the linear systems.
There are several potential explanations for this, including the fact that the
convergence radius of the Taylor expansion is limited, or that the convergence
of infGMRES is affected by the eigenvalues around.
However, increasing the iterations quadratically increases the cost of the
orthogonalization of the Arnoldi process and the cost of solving Hessenberg
systems, making the algorithm extremely slow.
Thus, it is too expensive to request highly accurate solutions by
performing many iterations of infGMRES ---
a wiser choice is to fix the infGMRES iterations and increase the
number of expansion points.
\end{remark}

\section{Numerical experiments}
\label{sec:numerexp}
In this section we present experimental results of
Algorithms~\ref{alg:wtinfgmres}.
All numerical experiments were performed using MATLAB R2022b on a Linux server
with two 16-core Intel Xeon Gold 6226R 2.90 GHz CPUs and 1024~GB main memory.

\subsection{Experiment settings}

\begin{figure}[tb!]
\centering
\includegraphics[height=10cm]{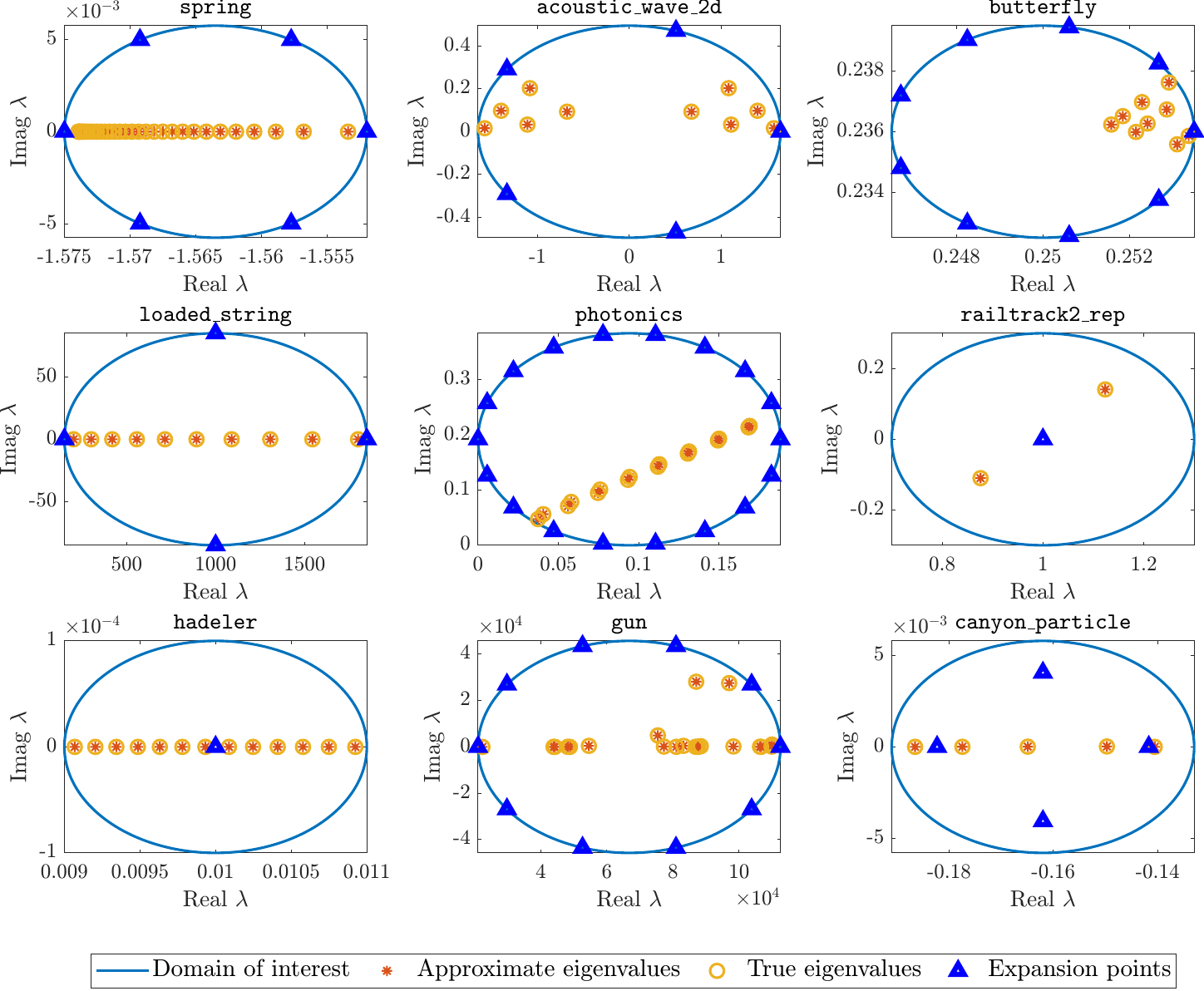}
\caption{Interest field, expansion points and eigenvalues of different test problems.}
\label{fig:eigpattern}
\end{figure}

Most of our test examples are chosen from the NLEVP collections~\cite{BHM2013}
except {\tt{photonics}}, which is similar to the one described
in~\cite{DAG2020}, but with a more general model for the permittivity as
in~\cite{GMG2017}.
We set the maximum number of the outer iterations of infGMRES as \(32\) for
all these examples, while the number of quadrature nodes and expansion points
varies from case to case.
Details regarding these examples, including their respective types, the
problem size \(n\), the number of quadrature nodes \(N\), and the number of
expansion points \(\nep\), are listed in Table~\ref{tab:expsinfo}.
For the exact distribution of these points or the contour, readers can check
Figure~\ref{fig:eigpattern}.

For computed approximate nonlinear eigenpair \((\hat\lambda,\hat v)\),
the convergence criterion is
\begin{equation}
\lVert T(\hat\lambda)\hat v\rVert_2\le
{\tt{tol}}\cdot\lVert T(\hat\lambda)\rVert_2\lVert\hat v\rVert_2,
\end{equation}
where {\tt{tol}} is the user-specified tolerance.
Unless otherwise stated, all the test results presented in this section
achieved an accuracy of \({\tt{tol}}= 10^{-12}\).

As we mentioned in Section~\ref{sub-sec:infGMRES}, it is possible to apply
the action \(T_0^{-1}\) in Algorithm~\ref{alg:wtinfgmres} in an inexact way,
e.g., using algebraic multigrid or GMRES, to make the algorithm even faster.
However, to emphasize the effect of our algorithm,
we use LU decomposition to solve all the linear systems
exactly in our numerical experiments.
Readers should keep in mind that, in practice,
inexact linear solvers can be implemented here for further acceleration.

\begin{table}[tb!]
\centering
\caption{Information of test problems. Here, \(n\) is the size of the problem
and \(k\) is the number of the eigenvalues to be computed.
The number of quadrature nodes, expansion points are represented by \(N\) and
\(\nep\), respectively.
The times consumed by Beyn's method with MATLAB backslash, \(t_{\rm{s}}\),
and with infGMRES \(t_{\rm{iG}}\) are also listed in the last columns.}
\begin{tabular}{c@{\hspace{0.8em}}c@{\hspace{0.8em}}c@{\hspace{0.8em}}c@{\hspace{0.8em}}c@{\hspace{0.8em}}c@{\hspace{0.8em}}c@{\hspace{0.8em}}c@{\hspace{0.8em}}c@{\hspace{0.8em}}c@{\hspace{0.8em}}c}
\hline
Problem & Type & \(n\) & \(k\) & \(N\) & \(\nep\) & \(t_{\rm{s}}\)(s) & \(t_{\rm{iG}}\)(s) & \((t_{\rm{s}}-t_{\rm{iG}})/t_{\rm{s}}\)\\
\hline
{\tt{spring}}             & QEP  & \(3000 \) & \(32\) & \(1024\) & \(6 \) & \(3.472\) & \(15.72\) & \(-353\%\)  \\
{\tt{acoustic\_wave\_2d}} & QEP  & \(9900 \) & \(10\) & \(512 \) & \(5 \) & \(28.14\) & \(9.475\) & \(66\%\)  \\
{\tt{butterfly}}          & PEP  & \(5000 \) & \(9 \) & \(512 \) & \(9 \) & \(11.54\) & \(8.453\) & \(27\%\)  \\
{\tt{loaded\_string}}     & REP  & \(20000\) & \(10\) & \(128 \) & \(4 \) & \(1.07 \) & \(9.805\) & \(-816\%\)  \\
{\tt{photonics}}          & REP  & \(20363\) & \(16\) & \(3060\) & \(18\) & \(600  \) & \(209.7\) & \(65\%\)  \\
{\tt{railtrack2\_rep}}    & REP  & \(35955\) & \(2 \) & \(128 \) & \(1 \) & \(533.4\) & \(71.95\) & \(87\%\)  \\
{\tt{hadeler}}            & NEP  & \(5000 \) & \(13\) & \(32  \) & \(1 \) & \(40.2 \) & \(46.6 \) & \(-16\%\)  \\
{\tt{gun}}                & NEP  & \(9956 \) & \(21\) & \(1024\) & \(10\) & \(501.2\) & \(148.9\) & \(70\%\)  \\
{\tt{canyon\_particle}}   & NEP  & \(16281\) & \(5 \) & \(256 \) & \(4 \) & \(40.94\) & \(7.652\) & \(81\%\)  \\
\hline
\end{tabular}
\label{tab:expsinfo}
\end{table}

\subsection{Overall performance on the test examples}

To illustrate the efficiency of our algorithm, we solve the linear systems by
MATLAB backslash (\(\backslash\)) as comparative experiments.
The times consumed by Beyn's method with MATLAB backslash (\(t_{\rm{s}}\))
and with infGMRES (\(t_{\rm{iG}}\)) are listed, respectively, in the last
columns of Table~\ref{tab:expsinfo}.
Except for the \({\tt{spring}}\), \({\tt{loaded\_string}}\) and
\({\tt{hadeler}}\), our algorithm achieved a speedup of at least \(27\%\).
The acceleration rate increases when the problems size becomes larger,
where our algorithm can benefit from taking fewer matrix decompositions
and reach a speedup over \(80\%\).

However, there are two scenarios where the advantages of our algorithm are
less evident.
In the cases of {\tt{spring}} and {\tt{loaded\_string}},
MATLAB implements specialized optimizations for tridiagonal matrices,
making it several times faster than the general MATLAB backslash operation.
On the other hand, in the case of {\tt{hadeler}},
two out of its three component matrices are dense.
Generally, iterative methods are not as efficient for dense matrices,
leading to reduced effectiveness in these specific instances.

\begin{figure}[tb!]
\centering
\includegraphics[height=4.6cm]{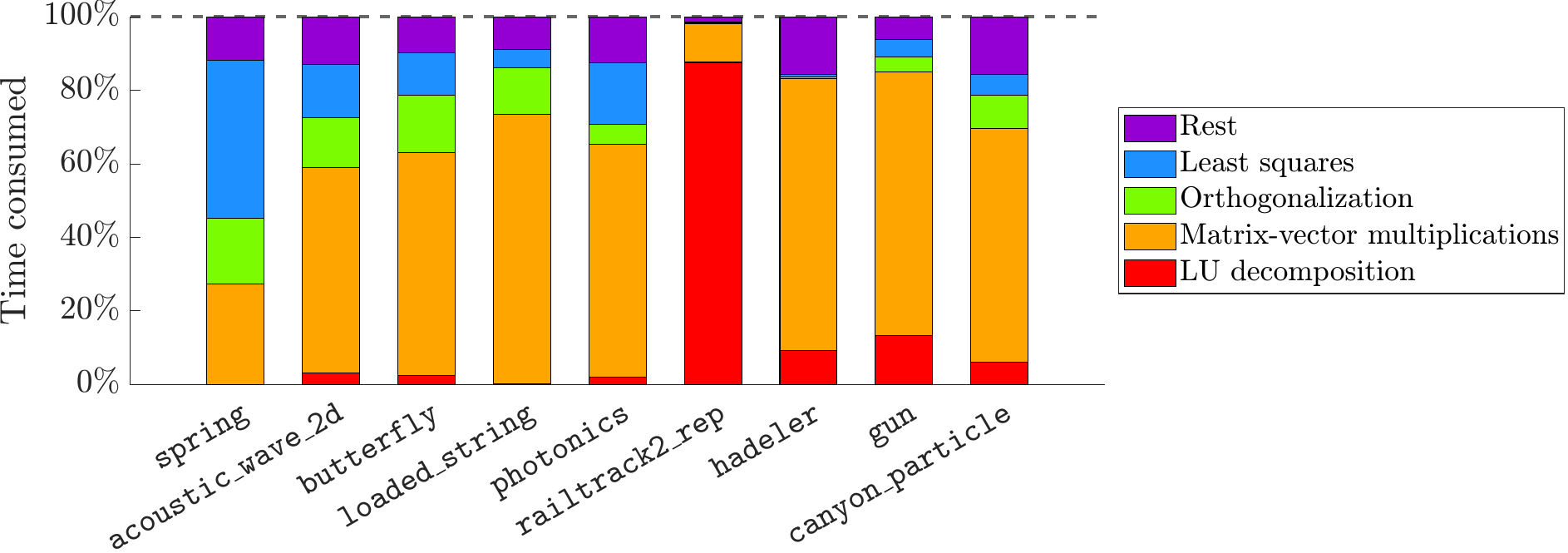}
\caption{The proportional chart of Beyn's methods with infinite GMRES.
The times consumed by different operations are illustrated.}
\label{fig:bar}
\end{figure}

To obtain a more detailed analysis, we have extracted the time proportions of
each operation, as illustrated in Figure~\ref{fig:bar}.
It can be found that a significant portion of time is dedicated to the
Arnoldi process (matrix--vector multiplications and orthogonalization) in our algorithm.
Even in the challenging problem {\tt{photonics}}, where over \(3000\)
quadrature nodes are needed to be solved, solving least squares problems
consume comparatively less time.
This implies that our algorithm proves especially advantageous in cases
involving a large number of quadrature nodes.

One thing to note from Figure~\ref{fig:bar} is that the time consumed by LU
decomposition is relatively higher for {\tt{railtrack2\_rep}}.
That is because the size of this problem is obviously larger, bringing about
more expensive factorizations.
However, it is actually a good news for our algorithm, which means that when
the problem size becomes larger, our algorithm can benefit more from solving
these factorizations in parallel.

\subsection{Comparison on different weighting strategies}

\begin{figure}[tb!]
\centering
\includegraphics[height=8cm]{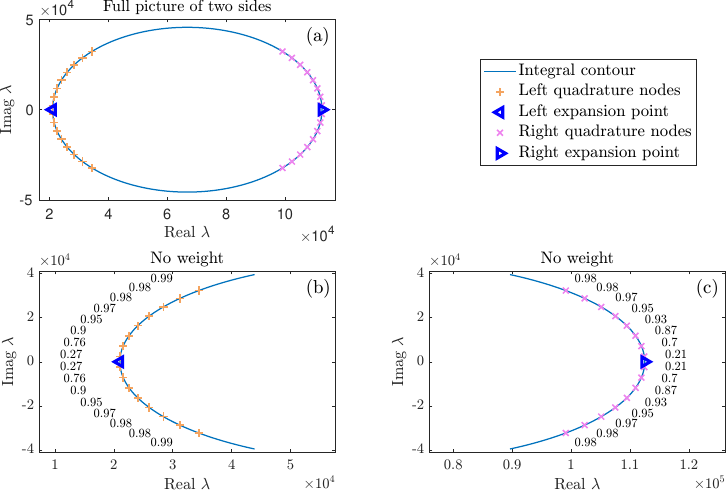}\\
\includegraphics[height=8cm]{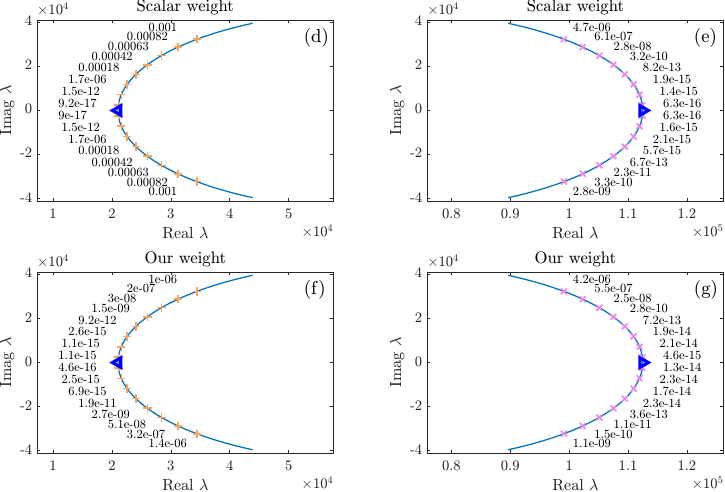}
\caption{Solving linear systems of the {\tt{gun}} problem
by infGMRES without weighting (b, c),
with scaling weighting~\eqref{eq:scalB2009}~(d, e)
and with our weighting~\eqref{eq:ourwt}~(f, g).
(a)~shows a full picture of the distribution of the expansion points and
part of the quadrature nodes.
The residuals plotted in the figures are the relative residuals of
solving linear systems.
For a certain quadrature node \(\xi_j\), the relative residual is defined as
\(\lVert T(\xi_j)x_{0,j}-z\rVert_2/\bigl(\lVert T(\xi_j)\rVert_2
\lVert x_{0,j}\rVert_2+\lVert z\rVert_2\bigr)\), where \(x_{0,j}\) stands for the approximate solution.}
\label{fig:compbalan}
\end{figure}

In Section~\ref{sec:weight}, we introduce a weighting strategy for infGMRES.
Here, we illustrate its efficiency under different circumstances with
numerical experiments.
The {\tt{gun}} problem is taken as the test problem.
We compute the Taylor expansion on the leftmost and rightmost points of the
contour, respectively, to approximate the linear systems corresponding to
quadrature nodes around.
The difference between these two points is that there are singularities close
to the left side one, leading to a more ill-conditioned problem.
The full picture can be found in Figure~\ref{fig:compbalan}~(a),
where we illustrate the corresponding relationships between expansion points
and quadrature nodes.
The number of outer iterations of infGMRES is fixed to \(32\) in all
experiments.
The residuals obtained with or without weighting can be found in
Figures~\ref{fig:compbalan}~(b), (c), (f), (g).
Our weighting strategy can bring significant boost on the accuracy in both cases.

As a comparison, we refer to the scaling strategy from~\cite[Section 6]{Betcke2009}.
In that algorithm, a scalar
\begin{equation}
\label{eq:scalB2009}
\rho=\left(\frac{\lVert T_0\rVert_2}{\lVert T_p\rVert_2}\right)^{1/p}
\end{equation}
is used to make the norm of \(\rho^j T_j\) as similar as possible.
Their consideration is that solving a PEP by applying a backward stable
algorithm is backward stable if
\[
\lVert T_0\rVert_2=\lVert T_1\rVert_2=\cdots=\lVert T_p\rVert_2.
\]
Even though this algorithm is designed for eigenvalue problems but not for
solving linear systems, it shares similar motivations with our algorithm.
We implement it in Figure~\ref{fig:compbalan}~(d), (e).
It can be found that, both strategies perform well for well-conditioned cases.
Nevertheless, the accuracy of~\eqref{eq:scalB2009} decays rapidly
in the ill-conditioned case; see Figure~\ref{fig:compbalan}~(d);
while our algorithm performs obviously better; see Figure~\ref{fig:compbalan}~(f).

\section{Conclusion}

In this work, we introduce infGMRES to reduce the cost of solving linear
systems in contour integral-based nonlinear eigensolvers.
We have worked out implementation details including the
convergence-accelerating weighting strategy, the memory-friendly TOAR-like
trick, and the selection of the parameters.
With these ingredients, we proposed a robust and efficient implementation of
infGMRES.
While our numerical experiments are carried out in Beyn's method, this
technique can actually be applied to all contour integral-based nonlinear
eigensolvers, where several moments are needed to be approximated.

Our method reduces the computational cost by reducing the number of required
matrix factorizations.
More precisely, it requires as many factorizations as the number of expansion
points, which is usually much smaller than the number of quadrature nodes.
This is especially relevant for difficult problems, where the quadrature rule
demands a large number of quadrature nodes or the scale is extremely large,
making a matrix decomposition very expensive.

Our future work may involve developing a machine learning-based adaptive
strategy for selecting expansion points automatically.
Additionally, we may implement a block variant of infGMRES to efficiently
solve all right-hand sides.

\section*{Acknowledgments}
We would like to thank Guillaume Dem\'{e}sy for providing us with the data from the photonics test case.
Furthermore, we would like to express our gratitude to the anonymous
reviewers.
Their feedback has been instrumental in enhancing this work.

Y.~Liu and M.~Shao were partly supported by the National Natural Science Foundation of China under grant No.~92370105.
J.~E.~Roman was supported by grant PID2022-139568NB-I00 funded by MICIU/AEI/10.13039/501100011033 and by ERDF/EU, and by grant RED2022-134176-T.
This work was carried out while Y.~Liu was visiting Universitat Polit\`ecnica de Val\`encia.

\end{document}